\documentclass[12pt]{article}
= 0.0in \headheight = 0.0in \topmargin = 0.3in
\def\Dj{\hbox{D\kern-.73em\raise.30ex\hbox{-}
\raise-.30ex\hbox{}}}
\def\dj{\hbox{d\kern-.33em\raise.80ex\hbox{-}
\raise-.80ex\hbox{\kern-.40em}}}
\usepackage{epsfig}
\usepackage{amsmath,amsthm,amsfonts,amssymb,amscd,cite}
\allowdisplaybreaks
\usepackage{mathrsfs}
\def \la{\lambda}

\newtheorem{theorem}{Theorem}

\newtheorem{lemma}[theorem]{Lemma}

\newtheorem{coro}[theorem]{Corollary}

\newtheorem{defi}[theorem]{Definition}

\newtheorem{remark}[theorem]{Remark}

\begin{document}

%

%
%

\title{On  bipartite graphs having minimal fourth adjacency coefficient\thanks{ Supported by Zhejiang Provincial Natural Science Foundation of
China (No. LY20A010005), and National Natural Science Foundation of
China (No. 11571315,11901525,11601006).}}
\author{Shi-Cai Gong\thanks{Corresponding author. E-mail addresses:
scgong@zafu.edu.cn(S. Gong).}~~and Shao-Wei Sun
\\{\small \it  School of Science, Zhejiang University of Science and Technology, }\\{\small \it
Hangzhou, 310023, P. R. China}
   }
\date{}
\maketitle

\baselineskip=0.20in

\noindent {\bf Abstract. } Let $G$ be a simple graph with order $n$ and adjacency matrix $\mathbf{A}(G)$. Let $\phi(G; \la)=\det(\la I-\mathbf{A}(G))=\sum_{i=0}^n\mathbf{a}_i(G)\la^{n-i}$ be the characteristic polynomial of $G$, where
 $\mathbf{a}_i(G)$ is called the $i$-th adjacency coefficient of $G$.
 Denote by $\mathfrak{B}_{n,m}$ the set of all connected graphs having $n$ vertices and $m$ edges.
A  bipartite graph  $G$ is referred as  bipartite optimal  if $$\mathbf{a}_4(G)=min\{\mathbf{a}_4(H)|H\in \mathfrak{B}_{n,m}\}.$$ The  value  $min\{\mathbf{a}_4(H)|H\in \mathfrak{B}_{n,m}\}$ is called the minimal $4$-Sachs number in $\mathfrak{B}_{n,m}$, denoted by $\bar{\mathbf{a}}_4(\mathfrak{B}_{n,m})$. \vspace{2mm}

For any given integer  pair $(n,m)$,  we in this paper investigate the bipartite optimal graphs. Firstly, we show that each bipartite optimal graph is a difference graph (see Theorem \ref{15}).    Then we deduce some structural properties on bipartite optimal graphs.
 As applications of those properties, we determine all bipartite optimal $(n,m)$-graphs together with the corresponding minimal $4$-Sachs number for $n\ge 5$ and $n-1\le m\le 3(n-3)$.
Finally, we  express the problem of computing the minimal $4$-Sachs number  as a class of combinatorial optimization problem, which relates to the partitions of positive integers.\vspace{3mm}



\baselineskip=0.30in

\section{Introduction}

Throughout the paper all  graphs are undirected and simple.  Let $G=(V,E)$ be a graph with vertex set $V=\{v_1,v_2,\ldots,v_n\}$ and edge set
$E$. The \emph{adjacency matrix} $\mathbf{A}=\mathbf{A}(G) = (a_{ij})_{n\times n}$ of $G$ is defined as
$a_{ij} = 1$ if and only if $v_i$ is adjacent to $v_j$, and $a_{ij} =0$ otherwise. 
 The \emph{ characteristic
polynomial} of  $G$, denoted by
$\phi(G;\la)$, is defined by $$\phi(G;\la)=\sum_{i=0}^n\mathbf{a_i}(G)\la^{n-i}=det(\la I-\mathbf{A}).$$
Hereafter, $\mathbf{a_i}(G)$ is  called  the $i$-th adjacency coefficient  of $G$.
  More results concerning the characteristic polynomial of graphs can be found in the literature \cite{big,br,cds,crs,gr,sw,m} and references therein. \vspace*{2mm}

 Let $G$ be a  graph. The subgraph  $H$ of $G$ is called a \emph{$p$-Sachs subgraph}  if the order of $H$ is $p$ and each component of $H$ is either a single edge or a cycle. Denote by $o(G)$ and $c(G)$ respectively the number of components and cycles contained in $G$. \vspace*{2mm}
For a general, not necessary be simple, graph $G$, as we known that the coefficient $\mathbf{a_i}(G)$  has a combinatorial interpretation
 in terms of $i$-Sachs subgraphs as follows; see example \cite[Theorem 1.2]{cds} $$\mathbf{a_i}(G)=\sum_H(-1)^{o(H)}2^{c(H)},\eqno{(1.1)}$$ where the
summation is over all $i$-Sachs subgraphs contained in $G$. Therefore, we sometimes refer $\mathbf{a_i}(G)$ as $i$-Sachs number for convenience.  \vspace*{2mm}

A graph having $n$ vertices and $m$ edges is referred as an $(n,m)$-graph. Denoted by $\mathfrak{B}_{n,m}$ the set of all connected bipartite $(n,m)$-graphs.
 An \emph{$r$-matching} in a given graph $G$ is
a subset with $r$ edges such that every vertex of $V(G)$ is incident
with at most one edge in it.
The  $r$-\emph{matching number}, denoted by $\mathbf{m}_r(G)$, is defined as the number of
$r$-matchings contained in $G$.\vspace*{2mm}

Let $G$ be a given graph. The $i$-th adjacency coefficients $\mathbf{a_i}(G)$   contains  abundant structural information and spectral information of such a graph obviously. Thus there have close relationships among them. For instance, if $G$ is acyclic, $\mathbf{m}_i(G)=(-1)^i\mathbf{a}_{2i}(G)$ for each $i(1\le i\le \lfloor\frac{n}{2}\rfloor)$;  if $G$ is bipartite, then $\mathbf{a_i}(G)=0$ for each odd number $i$ and $(-1)^{\frac{j}{2}}\mathbf{a}_j(G)\ge 0$ for each even number $j$; see for example \cite{big,cds,crs,lp}.\vspace*{2mm}

For any given $(n,m)$-graph $G$, from Eq. (1.1), we have $$\mathbf{a}_0(G)=1,~\mathbf{a}_1(G)=0,~\mathbf{a}_2(G)=-m,~\mathbf{a}_3(G)=-2c_3~~and~~\mathbf{a}_4(G)=\mathbf{m}_2(G)-2c_4,\eqno{(1.2)}$$
where $c_3$ and $c_4$ denote respectively the number of triangles and quadrangles contained in $G$. From (1.2), adjacency coefficients $\mathbf{a}_0(G)$, $\mathbf{a}_1(G)$ and $\mathbf{a}_2(G)$ are fixed, independent to the  structure of such a graph. Then it is interesting to investigate the relationship between the $i$-th adjacency coefficient $\mathbf{a}_i$ and the structural properties of a given graph. Moreover, $\mathbf{a}_3(G)=0$ if and only if $G$ contains no triangles. For more results concerning extremal triangle-free graphs, one can see  \cite{as,be,k} and references therein.
Naturally, it is interesting to study the relationships between $4$-Sachs number of a given graph and its structural properties.\vspace*{2mm}

A  bipartite graph  $G$ is referred as  bipartite optimal  if $$\mathbf{a}_4(G)=min\{\mathbf{a}_4(H)|H\in \mathfrak{B}_{n,m}\}.$$ The  value  $min\{\mathbf{a}_4(H)|H\in \mathfrak{B}_{n,m}\}$ is called the minimal $4$-Sachs number in $\mathfrak{B}_{n,m}$, denoted by $\bar{\mathbf{a}}_4(\mathfrak{B}_{n,m})$. \vspace{2mm}

In this paper, we will investigate the bipartite optimal graphs  and the corresponding minimal $4$-Sachs number $\bar{\mathbf{a}}_4(\mathfrak{B}_{n,m})$. The rest paper is organized as follows. In section $2$, we give some preliminary results, including the notation  threshold graphs,  difference graphs together with  their properties, and some other lemmas. In section $3$, we first give a compression operation that make graphs  minimize its $4$-Sachs number. Then we  show that each bipartite optimal graph is a difference graph. In section $4$,   we deduce some structural properties on bipartite optimal graphs.  As applications of those properties, we determine all bipartite optimal $(n,m)$-graphs together with the corresponding minimal $4$-Sachs number, for $n\ge 5$ and $n-1\le m\le 3(n-3)$, in section $5$. Finally,
 we experss the problem of computing the minimal $4$-Sachs number as a class of combinatorial optimization problem, which relates to the partitions of positive integers in section $6$.

%

%

\section{Preliminary}
Firstly, we introduce some preliminary results.
 Let $G=(V,E)$ be a  graph with $u\in V$.  Denote by $N_G(u)$ and $d_G(u)$ the \emph{neighbors} and the \emph{degree} of the vertex $u$, respectively. Vertices $u$ and $v$ of $G$ are called \emph{duplicate} if $N_G(u)=N_G(v)$. Let $G_1=(V_1,E_1)$ and
  $G_2=(V_2,E_2)$ be two graphs. Then
the  \emph{union} of $G_1$ and $G_2$, denoted by $G_1 \cup
G_2$, is defined as $(V_1 \cup V_2, E_1\cup E_2)$, and the \emph{join} of  $G_1$ and $G_2$, denoted by $G_1 \cap
G_2$, is defined as $(V_1 \cap V_2, E_1\cap E_2)$.
  Denote by $dis(u,v)$ the \emph{distance} between vertices $u$ and $v$.
Let $V_1\subset V$. The subgraph induced by the vertex set $V_1$ is denoted by $G[V_1]$. The cycle and the path of order $n$ are denoted by $C_n$ and $P_n$, respectively. The complete bipartite graph with bipartition $(X;Y)$ is denoted by  $K_{|X|,|Y|}$. The  graph  $K_{|X|,|Y|}$ is sometimes called a star if $min\{|X|,|Y|\}=1$.
\vspace{2mm}

\begin{defi} \label{01} {\em \cite{mp}} A graph $G = (V,E)$ is said to be a \emph{threshold} graph if there exists a threshold $t$ and a function $w : V(G) \rightarrow R$ such that $uv\in E(G)$  if and only if $w(u) +w(v) \ge t$.
\end{defi}
Threshold graphs have a beautiful structure and possess many important
mathematical properties such as being the extreme cases of certain graph
properties, see e.g.  \cite{mp,kr,lmt}. They also have applications in many areas such as computer
science and psychology. For more information on threshold graphs, one can see the book \cite{mp} and the references
therein.

\begin{defi} \label{02} {\em \cite{mp}} A graph $G = (V,E)$ is said to be  \emph{difference}  if there exists a threshold $t$ and a function $w : V(G) \rightarrow R$ such that $|w(v)| < t$ for all $v \in V$ and distinct vertices $u$ and $v$ are adjacent if and only if $|w(u) - w(v)| \ge t$.
\end{defi}

Difference graphs are called threshold bipartite graphs in \cite{mp} and chain graphs in \cite{y}. A threshold graph can be obtained from a difference graph by adding
all possible edges in one of the partite sets (on either side).
The following lemmas  are  useful to us.

\begin{lemma} \label{03} {\em \cite{kr}} The graph $G$ is difference if and only if $G$ is bipartite and the neighborhoods of vertices in one of the partite sets can be linearly ordered by inclusion.
\end{lemma}

\begin{lemma} \label{04} {\em \cite[Proposition 2.5(2)]{hps}} The connected bipartite graph $G$ is difference if and only if $G$  has no induced subgraph
$P_5$.
\end{lemma}

Let $G$ be a  difference graph with bipartition $(X; Y )$. Suppose that $X=\cup_{i=1}^kX_i$ and $Y=\cup_{i=1}^pY_i$ such that, for each $i$, both $X_i$ and $Y_i$ are non-empty, and all elements in $X_i$ (resp. $Y_i$) are duplicate. By Lemma \ref{03} we can further suppose that
  $$N(X_1)\supset N(X_2)\supset \ldots \supset N(X_k)~~ and ~~ N(Y_1)\supset N(Y_2)\supset \ldots \supset N(Y_p).$$ 
 Obviously, $p=k$. Furthermore, applying Lemma \ref{03}, $N_G(x)=\cup_{j=1}^{k-i+1}Y_j$ for  $x\in X_i$ and $N_G(y)=\cup_{j=1}^{k-i+1}X_j$ for  $y\in Y_i$. Thus
  for  each $i$ both $G[(\cup_{j=1}^{k-i+1}X_j)\cup Y_{i}]$ and $G[(\cup_{j=1}^{k-i+1}Y_j)\cup X_{i}]$ are complete bipartite.
Consequently,  the vertex set sequence $(X_1,X_2,\ldots,X_k; Y_1,Y_2,\ldots,Y_k)$ determines the difference graph $G$ and vice versa. Let $|X_i|=x_i$ and $|Y_i|=y_i$ for $i=1,2,\ldots,k.$
 For convenience, we refer $(X_1,X_2,\ldots,X_k; Y_1,Y_2,\ldots,Y_k)$ and $(x_1,x_2,\ldots,x_k; y_1,y_2,\ldots,y_k)$ as the \emph{vertex bipartition} and the \emph{vertex-eigenvector} of the difference graph $G$, respectively. The integer  $k$ above is called   the \emph{character} of  $G$.  Then  the complete bipartite graph $K_{n,m}$ is a difference graph with character $1$ and vertex-eigenvector $(n;m)$.\vspace{2mm}

Let $G$ be a  difference graph with vertex bipartition $(X_1,X_2,\ldots,X_k; Y_1,Y_2,\ldots,Y_k)$. Denoted by $G[\overline{(\cup_{j=1}^iX_j)\cup Y_{k-i+1}}]$
the subgraph  induced by the vertex set $V(G)\backslash (\cup_{j=1}^iX_j\cup Y_{k-i+1})$, named as the \emph{difference complement} of $G[(\cup_{j=1}^iX_j)\cup Y_{k-i+1}]$.
For each $i$, one can verify that $G[\overline{(\cup_{j=1}^iX_j)\cup Y_{k-i+1}}]$ is  a difference graph with vertex bipartition $(X_{i+1},\ldots,X_k;Y_1,\ldots,Y_{k-i})$.\vspace{2mm}

In the final of this section, we  give a preliminary Lemma, which will be used in Section $5$.\vspace{2mm}
\begin{lemma}\label{320} Let integers $n,m$ and $y$ satisfy $n>6$, $2(n-2)<m< 3(n-3)$ and $2<y<\frac{m+1-n}{2}.$ Then  $$\frac{(1+y)(n-1-y)-m}{y}>\frac{3n-9-m}{2}.$$
\end{lemma}
\noindent {\bf Proof.} Since $y>0$, it is sufficient to show that $$2(1+y)(y+1-n)+2m+y(3n-9-m)<0$$ holds for $2<y<\frac{m+1-n}{2}$, that is,  $$2y^2+(n-m-5)y +(2m-2n+2)<0$$ holds for $2<y<\frac{m+1-n}{2}$.
Let $$f(y)=2y^2+(n-m-5)y +(2m-2n+2).$$
Note that $f(2)=f(\frac{m+1-n}{2})=0$, then the result follows.\hfill $\blacksquare$
\section{ An operation}
 In this section, we first give a compression operation that make graphs  minimize their $4$-Sachs numbers. Then we show  that each  bipartite optimal graph  is a difference graph.\vspace{2mm}

  Let $u$ and $v$ be two vertices of the graph $G$. Define
 $$\mathbf{N}_G(u,v) = \{x\in V (G) \backslash \{u,v\} : xu\in E(G),xv\in E(G)\}$$ and
$$\mathbf{N}_G(u,\bar{v}) = \{x\in V (G) \backslash \{u,v\} : xu\in E(G),xv\notin E(G)\}.$$
Let $\mathbf{G}_{u\rightarrow v}$ be the graph formed by deleting all edges between $u$ and $\mathbf{N}_G(u,\bar{v})$ and adding all edges from $v$ to $\mathbf{N}_G(u,\bar{v})$. This operation is called the \emph{compression} of $G$ from $u$ to $v$; see Definition 2.4 in \cite{kr}. It is clear that  $\mathbf{G}_{u\rightarrow v}$ has the same number of edges as $G$.

Due to Keough and  Radcliffe \cite{kr}, a result on comparing the number of $k$-matchings between $G$ and $\mathbf{G}_{u\rightarrow v}$ is given as follows.\vspace{2mm}
 \begin{lemma} \label{12} \em{\cite[Lemma 4.1]{kr}} For all graphs $G$ and all $u,v\in V(G)$$$\mathbf{m_k}(G)\ge \mathbf{m_k}(G_{u\rightarrow v}).$$

\end{lemma}
Applying the method parallel to the proof of Lemma \ref{12}; see Lemma 4.1 in \cite{kr}, we can obtain a more strengthen result on counting the number of $k$-matchings, $k\ge 2$, of a graph. Since the proof is similar to that of Lemma \ref{12}, we omit the detail.
 \begin{lemma} \label{13} Let $G$ be a  graph and $u,v \in V(G)$. Then for any $k(k\ge 2)$  $$\mathbf{m}_k(G)\ge \mathbf{m}_k(G_{u\rightarrow v})$$
 inequality holds if and only if $\mathbf{N}_G(\bar{u},v)\neq \emptyset$ and $\mathbf{N}_G(u,\bar{v})\neq \emptyset.$
\end{lemma}

Applying Lemma \ref{13}, we have
 \begin{theorem} \label{14} Let $G$ be a graph with $u,v\in V(G)$.
If $dis(u,v)=2$, then
 $$\mathbf{a}_4(G)\ge \mathbf{a}_4(G_{u\rightarrow v})$$
 inequality holds if $\mathbf{N}_G(\bar{u},v)\neq \emptyset$ and $\mathbf{N}_G(u,\bar{v})\neq \emptyset.$
\end{theorem}
\noindent {\bf Proof.} Let $H:=G_{u\rightarrow v}$. Denote by  $Q(G)$  the set of all  quadrangles of $G$  and  set $q(G)=|Q(G)|$. From (1.2) $\mathbf{a}_4(G)=\mathbf{m}_2(G)-2q(G)$, then applying Lemma \ref{13} it is sufficiency to prove that  $$q(H)\ge q(G).\eqno{(3.1)}$$

To prove (3.1), we   construct an injection from $Q(G)\backslash Q(H)$ to $ Q(H)\backslash Q(G )$ that preserves
the number of quadrangles.
Firstly, we define a replacement function $r : E(G) \mapsto E(H)$  by $$r(e)=\left \{ \begin{array}{ll}
va, & {\rm if \mbox{ }}  e=ua\mbox{ } {\rm with} \mbox{ } a\in N_G(u);\\
ub, & {\rm if \mbox{ }}  e=vb\mbox{ } {\rm with} \mbox{ } b\in N_G(u,v);\\
e, & {\rm otherwise. \mbox{ }}
\end{array}\right.$$ Given $e\in E(G)$, we claim that $r(e)$ is an edge in $H$. If $y\in N_G(u)$, then $r(uy)=vy\in  E(H)$; if $y\in N_G(u,v)$, then $r(vy)=uy\in  E(H)$ and $r(e)=e\in E(H) $ if $e\in E(G)\backslash (E_1\cup E_2)$, where $E_1=\{ux|x\in N_G(u)\}$ and $E_2=\{vx|x\in N_G(u,v)\}$.\vspace{2mm}

Now we define an injection $\phi: Q(G)\backslash Q(H)\mapsto Q(H)\backslash Q(G)$   by
$$\phi(C)=\{r(e):e\in C, C\in Q(G)\backslash Q(H) \},$$
where $C$ is an arbitrary $4$-cycle of $Q(G)\backslash Q(H)$. Then $C$ must
  contain an edge $uw$ with $w\in N_G(u,\bar{v})$ and  another edge   $ux$ with $x\in N_G(u)$, regardless $x\in N_G(u,\bar{v})$ or $x\in N_G(u,v)$, that is, $C=uwyxu$ with $y\in N_G(w,x)$. By the definition of $r(e)$, $r(uw)=vw$, $r(ux)=vx$ and $r(e)=e$ if $e\notin \{uw,uy\}$. Then $\phi(C)=vwyxv$ and thus $\phi(C)\in Q(H)\backslash Q(G)$.\vspace{2mm}

It remain to  show that $\phi$ has a left inverse. Consider $r': E(H)\rightarrow E(G)$ defined by
$$r'(e)=\left \{ \begin{array}{ll}
ua, & {\rm if \mbox{ }}  e=va\mbox{ } {\rm with} \mbox{ } a\in N_G(u);\\
vb, & {\rm if \mbox{ }}  e=ub\mbox{ } {\rm with} \mbox{ } b\in N_G(u,v);\\
e, & {\rm otherwise. \mbox{ }}
\end{array}\right.$$
Define $\phi' : Q(H)\backslash Q(G)\rightarrow Q(G)\backslash Q(H)$ by $\phi'(C) = \{r¡ä(e) : e \in C\}$. It is straightforward to check that $\phi'(\phi(C)) = C$. Thus $\phi$ has a left inverse and so $\phi$ is injective. Consequently, the result follows.\hfill $\blacksquare$\vspace{2mm}
\vspace{2mm}

 \begin{remark} \label{141} Let $G$ be a graph with $u,v\in V(G)$. Then by the same method the result $\mathbf{a}_4(G)\ge \mathbf{a}_4(G_{u\rightarrow v})$ is also true if
$dis(u,v)>2$. The restriction ensure that the resultant graph $G_{u\rightarrow v}$ is connected.
\end{remark}


Combining with Lemmas \ref{03}, \ref{04} and Theorem \ref{14}, we have
 \begin{theorem} \label{15} Each bipartite optimal graph is a difference graph.
\end{theorem}
\noindent{\bf Proof.} Let $G$ be a bipartite optimal graph. From Lemma \ref{04}, $G$ is difference if and only if $G$ is $P_5$-free. Assume that $G$ contains the induced subgraph $P_5$, then there exist vertices $u$ and $v$ such that $u$ and $v$ lie in the same partite and satisfying $$N_G(u)\nsupseteq N_G(v)~~and ~~N_G(v)\nsupseteq N_G(u).$$ Applying Theorem \ref{14} $\mathbf{a}_4(G)>\mathbf{a}_4(G_{u\rightarrow v})$, which is a contradiction to Lemma \ref{03}.\hfill $\blacksquare$\vspace{2mm}


\section{Computing the minimal $4$-Sachs number in $\mathfrak{B}_{n,m}$}

In this section, we study the problem of computing the  minimal $4$-Sachs number in $\mathfrak{B}_{n,m}$.
From Theorem \ref{15},  each  bipartite optimal graph is difference. Henceforth, we focus   on difference graphs. We begin our discussion with   a formula on  $4$-Sachs number  of a  difference graph.
\vspace{2mm}

Let  $G$ be a graph,  $C$  an even cycle of  $G$ with length
$2l$ and $H$   a $2r$-Sachs subgraph of $G$.
Suppose that $r\ge l$. We say the cycle $C$  is \emph{embedded} in
$H$ if  $C\cap H$ forms a $2l$-Sachs subgraph  and $C\cup H$ forms a $2r$-Sachs subgraph; see \cite{g}. Applying the formula (1.1), we have

   \begin{lemma}\label{31} Let $G$ be a bipartite graph and $C_4$  a given $4$-cycle of  $G$. Denote by $\mathbb{H}(C_4,2r)$ the set
of all  $2r$-Sachs subgraphs, of $G$, embedding the  cycle $C_4$. Then $$\sum_{H\in \mathbb{H}(C_4,2r)}(-1)^{o(H)}2^{c(H)}=0,$$ where the summation is over all $2r$-Sachs subgraphs of $\mathbb{H}(C_4,2r)$.
\end{lemma}
\noindent {\bf Proof.} Obviously, $r\ge 2$. Let $C_4=x_1y_1x_2y_2x_1$. Since $G$ is bipartite,  $G[\{x_1,y_1,x_2,y_2\}]=C_4$.  If  $r=2$, then $\mathbb{H}(C_4,2r)$ contain exactly three elements:
two disjoint $2$-matchings of $C_4$, named as $M_1=\{x_1y_1, x_2y_2\}$ and $M_2=\{y_1x_2, y_2x_1\}$,  and $C_4$ itself. Thus
$$\sum_{H\in \mathbb{H}(C_4,2r)}(-1)^{o(H)}2^{c(H)}=(-1)^{o(C_4)}2^{c(C_4)}+(-1)^{o(M_1)}2^{c(M_1)}+(-1)^{o(M_2)}2^{c(M_2)}=0.
$$
If $r>2$, then each $H\in \mathbb{H}(C_4,2r)$ contains either $M_1$ or $M_2$ or $C_4$ as a subgraph. Let $H=H_1\cup H_2$, where $H_1$ is the $(2r-4)$-Sachs subgraph of $G\backslash C_4$ and $H_2=\{M_1,M_2,C_4\}$. Thus
 \begin{equation*}
\begin{array}{lll} \sum_{H\in \mathbb{H}(C_4,2r)}(-1)^{o(H)}2^{c(H)}\\=\sum_{H_1\in G\backslash C_4}(-1)^{o(H_1)}2^{c(H_1)}[(-1)^{o(C_4)}2^{c(C_4)}+(-1)^{o(M_1)}2^{c(M_1)}+(-1)^{o(M_2)}2^{c(M_2)}]\\=0.
\end{array}
\end{equation*}
Consequently, the result follows.\hfill $\blacksquare$\vspace{2mm}

  As a consequence of Lemma \ref{31}, a  formula on  $4$-Sachs number  of  difference graphs can be obtained.\vspace{2mm}

   \begin{theorem}\label{32} Let $(X_1,X_2,\ldots,X_k; Y_1,Y_2,\ldots,Y_k)$$(k\ge 1)$ be the vertex bipartition of the difference graph $G$.  Then
   $$\mathbf{a}_{4}(G)=\sum_{i=1}^{k-1}\mathbf{a}_2(G[X_i;\cup_{j=1}^{k-i+1}Y_j])\mathbf{a}_{2}(\overline{G[X_i;\cup_{j=1}^{k-i+1}Y_j]}).$$
\end{theorem}
\noindent {\bf Proof.} Recall that  $\mathbf{a}_2(G)$ is the opposite of the number of edges contained in $G$ by Eq.(1.2). By the discussion above, $E(G)$ can be partitioned as
$$\bigcup_{i=1}^k E(G[X_i;\cup_{j=1}^{k-i+1}Y_j]).$$
Applying Lemma \ref{31}, to compute $\mathbf{a}_{4}(G)$, it is sufficiency to count the number of all $2$-matchings in which no two edges are contained in any quadrangle. Let $ M$ be such an $2$-matching. If $e\in G[X_i;\cup_{j=1}^{k-i+1}Y_j]$ $(i=1,2,\ldots,k-1)$, then the another edges of $H$ must contained  in $G\overline{[X_i;\cup_{j=1}^{k-i+1}Y_j]}$. Conversely, each pair edges $(e_1,e_2)$ with $e_1\in G[X_i;\cup_{j=1}^{k-i+1}Y_j]$ and $e_2\in G[\overline{X_i;\cup_{j=1}^{k-i+1}Y_j}]$ forms a $2$-matching  which does not embedded in  any quadrangle. Consequently, the result follows.\hfill $\blacksquare$\vspace{2mm}

Let $(X_1,X_2,\ldots,X_k; Y_1,Y_2,\ldots,Y_k)$$(k\ge 1)$ be the vertex bipartition of the difference graph $G$.  By the symmetry, we have
   $$\mathbf{a}_{4}(G)=\sum_{i=1}^{k-1}\mathbf{a}_2(G[\cup_{j=1}^{k-i+1}X_j;Y_i])\mathbf{a}_{2}(\overline{G[\cup_{j=1}^{k-i+1}X_j;Y_i]}).$$
   Therefore, without loss of generality, we in the following always suppose that
        $$ \sum_{i=1}^k|X_i|\ge \sum_{i=1}^k|Y_i|.$$\vspace{2mm}

    Difference graphs can be represented by Young diagrams intuitively \cite{kr}. 

\begin{defi} \label{33} Let $(X_1,X_2,\ldots,X_k; Y_1,Y_2,\ldots,Y_k)$ be the vertex-bipartition of the difference graph $G$. The Young diagram, or Young matrix, $Y(G)=(y_{ij})$ is defined as follows:
First, we set the rows of $Y$ correspond to the vertices $x^1_{1},\ldots,x^{|X_1|}_{1},\ldots,x^1_{k},\ldots ,x^{|X_k|}_{k}$ and the columns of $Y$ correspond to the vertices $y^1_{1},\ldots,y^{|Y_1|}_{1},\ldots,y^1_{k},\ldots ,y^{|Y_k|}_{k}$, respectively. Then we define $y_{ij}=1$ if and only if the vertices corresponding to the row $i$ and the column $j$ are adjacent, and $y_{ij}=0$ otherwise.
\end{defi}

 To compute the minimal  $4$-Sachs number, we need to introduce another matrix, named as  the characteristic matrix, as follows.

\begin{defi} \label{34} Let $(x_1,x_2,\ldots,x_k; y_1,y_2,\ldots,y_k)$ be the vertex-eigenvector of the  difference graph $G$. The characteristic matrix $T(G)=(t_{ij})_{k \times k}$ is defined as follows: $t_{ij}=x_iy_j$ if $i+j\le k+1$ and $t_{ij}=0$ otherwise, that is, $$T=\left(
      \begin{array}{ccccc}
        x_1y_1 & x_1y_2 & \cdots & x_1y_{k-1} & x_1y_k \\
         x_2y_1 & x_2y_2 & \cdots & x_2y_{k-1} & 0 \\
        \cdots & \cdots & \cdots & \cdots & \cdots \\
        x_{k-1}y_1 & x_{k-1}y_2& \cdots & 0 & 0 \\
         x_{k}y_1 & 0 & \cdots & 0 & 0 \\
      \end{array}
    \right).$$
\end{defi}\vspace{2mm}

 Let $A=(a_{ij})_{n\times n}$ be a matrix and $S,T$ be two sub-index sets of  $\{1,2,\ldots,n\}$. Set $\{1,2,\ldots,n\}=:\langle n \rangle$ and $\bar{S}:=\langle n \rangle \backslash S$. Denote by $A(S;T)$  the submatrix of $A$ by deleting the rows indicated by $\overline{S}$ and the columns indicated by $\overline{T}$. The column matrix $A(\langle n \rangle;\{i\})$ will be written as $A(\cdot ;i)$ for simplify. In addition, we use $s(A)$ to denote the sum of all entries of $A$. Then we have

  \begin{theorem} \label{37} Let $(x_1,x_2,\ldots,x_k; y_1,y_2,\ldots,y_k)$ be the vertex-eigenvector of the  difference graph $G$. Suppose that   the characteristic matrix of  $G$ is $T=(t_{ij})_{k\times k}$.  Then
  $$\mathbf{a}_4(G)=
\sum_{i=1}^{k-1} s(T(\cdot ; k-i+1))s(T(\overline{\langle i\rangle};\langle k-i\rangle))$$
\end{theorem}
\noindent{\bf Proof.}  Applying Theorem \ref{32} we need only  to count the number of those $2$-matchings in which each of them does not embedded in any quadrangle. Let $M$ be such a matching. If $e\in M$ is contained in $ G[\cup_{j=1}^iX_j;Y_{k-i+1}]$, then the another  edge of $M$ must contained  in $G[\overline{\cup_{j=1}^iX_j;Y_{k-i+1}}]$. Note that the number of edges contained in $G[\cup_{j=1}^iX_j;Y_{k-i+1}]$ is $s(T(\cdot;k-i+1))$ and the number of edges contained in $G[\overline{\cup_{j=1}^iX_j;Y_{k-i+1}}]$ is $s(T(\overline{\langle i\rangle};\langle k-i\rangle))$. Thus the result follows.\hfill $\blacksquare$\vspace{2mm}

Based on Theorem \ref{37}, we can deduce some properties on the vertex-eigenvector of the bipartition optimal graphs.

 \begin{theorem} \label{312} Let $G$ be a bipartition optimal graph  in $\mathfrak{B}_{n,m}$.
 Let also $(x_1,x_2,x_3,\ldots,x_k; \\ y_1,y_2,\ldots,y_k)$($k\ge 2$) be its vertex-eigenvector. Suppose that $\sum_{i=1}^{k}x_i\ge  \sum_{j=1}^{k-1}y_j.$ Then $$x_1> y_1.$$
\end{theorem}
\noindent {\bf Proof.} Assume to the contrary that $x_1\le y_1,$ say $y_1= x_1+y_1^*$ with $y_1^*\ge 0$.
 Let $G_1$ be the difference graph with vertex-eigenvector  $(x_1+x_2,x_3,\ldots,x_k,y_k; x_1,y_1^*,y_2,\ldots,y_{k-1})$.
Then
$$T_1:=T(G)=\left(
                        \begin{array}{ccccc}
                          x_1(x_1+y_1^*) & x_1y_2 & \cdots & x_1y_{k-1} & x_1y_k \\
        x_2(x_1+y_1^*) & x_2y_2 & \cdots & x_2y_{k-1} & 0 \\
        \cdots & \cdots & \cdots & \cdots & \cdots \\
        x_{k-1}(x_1+y_1^*) & x_{k-1}y_2& \cdots & 0 & 0 \\
         x_{k}(x_1+y_1^*) & 0 & \cdots & 0 & 0 \\

                        \end{array}
                      \right)$$ and
               $$T_2:=T(G_1)=\left(
                        \begin{array}{ccccc}
        (x_1+x_2)x_1& (x_1+x_2)y^*_1 & (x_1+x_2)y_2 & \cdots & (x_1+x_2)y_{k-1}  \\
       \cdots & \cdots & \cdots & \cdots & \cdots \\
        x_{k-1}x_1&x_{k-1}y^*_1 & x_{k-1}y_2& \cdots & 0 \\
        x_{k}x_1& x_{k}y^*_1 & 0 & \cdots & 0 \\
        x_1y_k& 0 & 0 & \cdots & 0  \\
                        \end{array}
                      \right).$$
One find that both $T_1$ and $T_2$ are square matrices with order $k$. Applying Theorem \ref{37},
 we have
\begin{equation*}
\begin{array}{lll}\mathbf{a}_4(G)-\mathbf{a}_4(G_1)&=&\sum_{i=1}^{k-1}s(T_1(\cdot; k-i+1))s(T_1(\overline{\langle i \rangle};\langle k-i \rangle))\\&&-\sum_{i=1}^{k-1}s(T_2(\cdot; k-i+1))s(T_2(\overline{\langle i \rangle};\langle k-i \rangle))
\\
&=&x^2_1y_k(\sum_{i=2}^{k}x_i-\sum_{j=2}^{k-1}y_j-y_1^*)\\
&=&x^2_1y_k(\sum_{i=1}^{k}x_i-\sum_{j=1}^{k-1}y_j)\\
&> &0,
\end{array}
\end{equation*}
 which yields a contradiction to the minimality of $4$-Sachs number  of $G$. Consequently, the result follows.\hfill $\blacksquare$\vspace{2mm}

Further, we have

 \begin{theorem} \label{3121} Let $G$ be a bipartition optimal graph  in $\mathfrak{B}_{n,m}$.
 Let also $(x_1,x_2,x_3,\ldots,x_k; \\ y_1,y_2,\ldots,y_k)$($k\ge 3$) be its vertex-eigenvector. Suppose that $\sum_{i=1}^{k}x_i\ge  \sum_{j=1}^{k-1}y_j,$ then $$x_1\ge y_1+y_2.$$
\end{theorem}
\noindent {\bf Proof.} By Theorem \ref{312} the result follows if $y_2=1$. Suppose now that $y_2\ge 2$. Assume to the contrary that $y_1<x_1< y_1+y_2,$ say $x_1= y_1+y_2^*$ and $y_2=y_2^*+ y_2^{**}$ with $y_2^*, y_2^{**}\ge 1$.
 Let $T_1:=T(G)$ and $G_1$ be the difference graph with vertex-eigenvector  $(x_1+x_2,x_3,\ldots,x_{k-1},y_k, x_k; y_1,y_2^*,y_2^{**},y_3,\ldots,y_{k-1})$.
Then one can verify that $G$ and $G_1$ have the same number of edges, and
               $$T_2:=T(G_1)=\left(
                        \begin{array}{ccccc}
       (x_1+x_2)y_1& (x_1+x_2)y^*_2 & (x_1+x_2)y^{**}_2 & \cdots & (x_1+x_2)y_{k-1}  \\
       \cdots & \cdots & \cdots & \cdots & \cdots \\
        x_{k-1}y_1&x_{k-1}y^*_2 & x_{k-1}y^{**}_2& \cdots & 0 \\
        y_{k}y_1& y_{k}y^*_2 & 0 & \cdots & 0 \\
        x_1y_1& 0 & 0 & \cdots & 0  \\
                        \end{array}
                      \right).$$
 Thus, we have
\begin{equation*}
\begin{array}{lll}\mathbf{a}_4(G)-\mathbf{a}_4(G_1)&=&x_1y_k[(y_1+y_2^*)\sum_{i=2}^kx_i-x_1\sum_{i=3}^{k-1}]-x_1^2y_ky_2^{**}-x_ky_1y_ky_2^*
\\
&=&x^2_1y_k(\sum_{i=2}^{k}x_i-\sum_{i=2}^{k-1}y_i)+x_1^2y_ky_2^*-x_ky_1y_ky_2^*\\
&=&x^2_1y_k(\sum_{i=1}^{k-1}x_i-\sum_{i=1}^{k-1}y_i)+x_1^2y_ky_2^*+x_1^2y_kx_k-x_ky_1y_ky_2^*\\
&> &0,
\end{array}
\end{equation*}
 which yields a contradiction to the minimality of $4$-Sachs number  of $G$. Consequently, the result follows.\hfill $\blacksquare$\vspace{2mm}

  In addition, we need to define a compression move that makes  difference graphs having more minimal $4$-Sachs number.\vspace{2mm}
 \begin{defi}\label{314} Let  $Y=(y_{ij})_{m\times n}$ be the Young matrix of the difference graph $G$. The entry $(i,j)$ is called out-corner if $y_{i,j}=1$ and $ y_{i+1,j}=y_{i,j+1}=0$.
 The entry $(p,q)$ is called an in-corner if $y_{p-1,q}=y_{p,q-1}=1$ and $ y_{p,q}=0$. If $(i,j)$ is a out-corner, $(p,q)$ is an in-corner, and $(i,j)$ and $(p,q)$ are not adjacent, we use $Y_{ij\rightarrow pq}$ to denote the matrix obtained from $Y$ by setting $y_{i,j}=0$ and $y_{p,q}=1$. This is called the difference compression of $Y$ from $(i,j)$ to $(p,q)$.
\end{defi}
It is clear that $Y_{ij\rightarrow pq}$ is also the Young matrix of a difference graph and those two  graphs have the same number of edges.

 \begin{lemma} \label{315} Let  $Y=(y_{ij})_{m\times n}$ be the Young matrix of the difference graph $G$. Suppose that $(i,j)$ is a out-corner vertex of $Y$. Denote by $e$ the edge of $G$ corresponding to the entry $(i,j)$.  Then   the number of $2$-matchings containing the edge $e$ and  embedding no $(even)$ cycles equals $$s(Y)-i j.$$
\end{lemma}
\noindent {\bf Proof.}   Let $e'$ be an edge of $G$ whose corresponding entry in $Y$ is $(p,q)$ such that $M=\{e',e\}$ is a  $2$-matching    embedding no $(even)$  cycles. Then by Lemma \ref{31} $p>i$ or $q>j$, that is, each entry of the submatrix $Y(\langle i\rangle;\langle i\rangle)$ is not contained. Thus the result follows.\hfill $\blacksquare$\vspace{2mm}
 \begin{theorem} \label{316} Let  $Y$ and $Y'$ be the Young matrices of the difference graphs $G$ and $G'$, respectively.  Suppose that $Y'=Y_{ij\rightarrow pq}$,  the  compression of $Y$ from $ij$ to $pq$. Then $$\mathbf{b}_4(G)\ge \mathbf{b}_4(G')$$ if and only if $ij\le pq$ and inequality holds if $ij< pq$.
\end{theorem}
\noindent {\bf Proof.} By Definition \ref{314}, $(i,j)$ is a out-corner and $(p,q)$ is an in-corner of $Y$. Denote by $Y^*$ the matrix obtained from $Y$ by replacing the entry $y_{pq}$ by $1$. One can verify that $Y^*$ is also the Young matrix of a difference graph, denoted by $G^*$. Then it is sufficiency to show that the cardinality of $M(ij)$ is no less than that of $M(pq)$, where  $M(ij)$ (resp. $M(pq)$) denotes all $2$-matchings containing the edge $ij$ (resp. $pq$) and  embedding no even cycles.\vspace{2mm}

By Lemma \ref{315}$$\mathbf{a}_4(G)- \mathbf{a}_4(G')=ij-pq.$$ Thus the result follows.\hfill $\blacksquare$\vspace{2mm}

By the method similar to Theorem \ref{316}, we have

 \begin{coro} \label{317} Let  $Y$ be the Young matrix of the difference graph $G$. Let $\{P_i(a_i,b_i)|i=1,2,\ldots,s\}$ and $\{Q_i(c_i,d_i)|i=1,2,\ldots,s\}$ be two vertex sequences. Let also $G_0=G$ and $G_i=G_{i-1}-P_i+Q_i$ for $i=1,2,\ldots,s$. Suppose that for each $i$ $P_i$ is a outer corner and $Q_i$ is an inner corner of $G_{i-1}$. Then $\mathbf{a}_4(G)> \mathbf{a}_4(G_s)$ if $$\sum_{i=1}^sa_ib_i>\sum_{i=1}^sc_id_i.$$
\end{coro}\vspace{2mm}

\section{Bipartite optimal graphs}
Applying all  preliminary results above, we determine some bipartite optimal graphs together with the corresponding minimal $4$-Sachs number.
First of all,  The following result is obviously.
 \begin{theorem} \label{311} Let positive integers $t$, $n$ and $m$ satisfy  $m=t(n-t)$. Then the complete bipartite graph $K_{t,n-t}$ is the unique  bipartite optimal graph in $\mathfrak{B}_{n,m}$.
\end{theorem}

Therefore, we focus on those integer pair $(n,m)$ satisfying $t(n-t)<m<(t+1)(n-t-1)$ for some integer $t$. Especially, we have
 \begin{theorem} \label{319} Let $n\ge 6$ and $n-1<m< 2(n-2)$. Then  the  difference graph with vertex eigenvector $(1,1;m-n-2,2n-4-m)$ is the unique bipartite optimal graph in $\mathfrak{B}_{n,m}$.
\end{theorem}
\noindent{\bf Proof.} Let $G$ be the   bipartite optimal graph in $\mathfrak{B}_{n,m}$.  By Theorem \ref{15}, $G$ is difference. Suppose that the vertex-eigenvector of $G$ is $(x_1,x_2,\ldots, x_k; y_1,y_2,\ldots, y_k )$. Since $n-1<m< 2(n-2)$, $k\ge 2$. By Theorem \ref{311} $y_1=1$. Then $$\mathbf{a}_4(G)\ge x_k(m-\sum_{i=1}^kx_i)$$ with equality  if and only if $k=2.$ Moreover, $\sum_{i=1}^{k-1}x_i\le \frac{m-x_k}{2}$ with equality  if and only if $\sum_{i=1}^ky_i=2,$ then $m-\sum_{i=1}^kx_i\ge \frac{m-x_k}{2}$ with equality  if and only if $\sum_{i=1}^ky_i=2.$ Consequently, $$\mathbf{a}_4(G)\ge \frac{x_k (m-x_k)}{2}=(2n-4-m)(m-n+1)$$ with equality  if and only if $k=2$ and $y_2=1.$ Thus the vertex eigenvector of $G$ is $(1,1;m-n-2,2n-4-m)$, whose character is $2$. Consequently, the result follows.\hfill $\blacksquare$\vspace{2mm}

 \begin{theorem} \label{320}Let $n>6$ and $2(n-2)<m< 3(n-3)$. Let  $G$ be  a bipartite optimal graph in $\mathfrak{B}(n,m)$. Then $G$ is a  difference graph and the corresponding vertex eigenvector $w$  satisfies
 $$w=\left \{ \begin{array}{ll}
(m-2n+6,3n-m-9;2,1), & {\rm if  \mbox{ } m< \frac{7n}{3}-7 };\\
(\frac{n-3}{3},\frac{2n-6}{3};2,1)~or~(\frac{2n-6}{3},\frac{n-3}{3};1,2), & {\rm if  \mbox{ } m=\frac{7n}{3}-7 };\\
(\frac{m-n+3}{2},\frac{3n-m-9}{2};1,2), & {\rm if  \mbox{ } m> \frac{7n}{3}-7 \mbox{ } and \mbox{ }3n-m-9\mbox{ }is \mbox{ } even};\\
(\frac{m-n+2}{2},1,\frac{3n-m-10}{2};1,1,1), & {\rm if  \mbox{ } m> \frac{7n}{3}-7 \mbox{ } and \mbox{ }3n-m-9\mbox{ }is \mbox{ } odd.} \end{array}\right.$$
\end{theorem}

\noindent{\bf Proof.} Before beginning our proof, we should point out that $3n-m-9$ is always even if $m=\frac{7n}{3}-7.$ In addition, we  sometimes use $\mathbf{a}_4(G,w)$ to denote the $4$-Sachs number of the difference graph $G$ with vertex-eigenvector $w$ for differentiation.\vspace{2mm}

  By Theorem \ref{15}, $G$ is difference. Suppose that the vertex-eigenvector of $G$ is  $(x_1,x_2\ldots, x_k;$ $ y_1,y_2\ldots, y_k )$. Without loss of generality, suppose that $\sum_{i=1}^kx_i\ge \sum_{i=1}^ky_i$. Recall that $2(n-2)<m< 3(n-3)$, then $y_1\le 2$ by Theorem \ref{311}.\vspace{2mm}

If $y_1=2$, then $$\mathbf{a}_4(G)\ge 2 x_k(m-2\sum_{i=1}^kx_i)$$ with equality   if and only if $k=2$. Moreover,  $\sum_{i=1}^{k-1}x_i\le \frac{m-2x_k}{3}$ with equality  if and only if $\sum_{i=1}^ky_i=3.$ Thus $m-2\sum_{i=1}^kx_i\ge \frac{m-2x_k}{3}$ with equality  if and only if $\sum_{i=1}^ky_i=3.$
 Consequently, $$\mathbf{a}_4(G)\ge \frac{2x_k (m-2x_k)}{3}=2(3n-9-m)(m-2n+6)$$ with equality  if and only if $k=2$ and $y_2=1.$ Then the   vertex-eigenvector of $G$  is $w_1=(m-2n+6,3n-9-m; 2,1)$, that is, $$\mathbf{a}_4(G,w_1)=2(3n-9-m)(m-2n+6).\eqno{(5.1)}$$

If $y_1=1$, then $$\mathbf{a}_4(G)\ge x_k(m-\sum_{i=1}^kx_i)$$ with equality   if and only if $k=2$. 
 Recall that $2(n-2)<m< 3(n-3)$,  then $\sum_{i=1}^ky_i\ge 3$ and  $\sum_{i=1}^{k-1}x_i\le \frac{m-x_k}{3}$ with equality   if and only if $\sum_{i=1}^ky_i= 3$. Consequently,
$$\mathbf{a}_4(G)\ge x_k(m-\sum_{i=1}^kx_i)=\frac{2x_k(m-x_k)}{3}=\frac{(3n-9-m)(m-n+3)}{2}$$ with equality  if and only if $k=2$ and $y_2=2$. In such a case the  vertex-eigenvector of $G$ is $w_2=(\frac{m-n+3}{2},\frac{3n-9-m}{2}; 1,2)$, that is, $$\mathbf{a}_4(G,w_2)=\frac{2x_k(m-x_k)}{3}=\frac{(3n-9-m)(m-n+3)}{2}.\eqno{(5.2)}$$

The condition that $(\frac{m-n+3}{2},\frac{3n-9-m}{2}; 1,2)$ being of the vertex-eigenvector of $G$   compels
 that $3n-m-9$ is even as $\frac{3n-9-m}{2}$ is integer. Then it remain to consider the case that $y_1=1$ and $3n-m-9$ is odd.\vspace{2mm}

  For $y_1=1$ and $3n-m-9$ is odd, we below  divide our proof into five assertions to show   that the vertex-eigenvector of the desired difference graph $G$  is $w_3=(\frac{m-n+2}{2},1,\frac{3n-m-10}{2};1,1,1)$, and the corresponding $4$-Sachs number is
$$\mathbf{a}_4(G, w_3)=
    \frac{(3n-m-10)(m-n+3)}{2}+m-n+2.\eqno{(5.3)}$$

\noindent {\bf Assertion 1.} The character  is no less than $3$.\\
Assume to the contrary that the character  is   $2$. Then $$\mathbf{a}_4(G)= x_1x_2y_2.$$
 Recall that $3n-m-9$ is odd, then $y_2>2$. Moreover,  if $m\le \frac{7n}{3}-7$, then $x_1y_2>m-n+3>4(m-n+6)$ and $x_2=\frac{(1+y_2)(n-1-y_2)-m}{y_2}\ge \frac{3n-9-m}{2}$ by Lemma \ref{320}. Thus $$\mathbf{a}_4(G)>2(3n-9-m)(m-2n+6)=\mathbf{a}_4(G,w_1),$$ which is a contradiction. If $m> \frac{7n}{3}-7$, then
$$\left \{ \begin{array}{ll}
x_1y_2= m-n+1+y_3\ge m-n+4;\\
x_2=n-2-y_2-\frac{m-n+1}{y_2}\ge \frac{4n-m-16}{3} \end{array}\right.$$
 with equality  if and only if $y_2=3$. Thus
 $$\mathbf{a}_4(G)\ge \frac{(m-n+4)(4n-m-16)}{3}$$with equality  if and only if $y_2=3$.
  Recall that $n>6$ and $m>\frac{7n}{3}-7$, then
\begin{equation*}
\begin{array}{lll}\mathbf{a}_4(G)-\mathbf{a}_4(G,w_3)= \frac{(m-n+3)(m-n-6)}{6}+n-\frac{10}{3}>0,
\end{array}
\end{equation*}
which is also a contradiction.

Therefore,  suppose that the  vertex-eigenvector of $G$ is $w=(x_1,x_2,\ldots,x_k;1,y_2,\ldots,y_k)$ with $k\ge 3$.

\noindent {\bf Assertion 2.} $y_2=1$.
Assume to the contrary that $y_2\ge 2$.  By Theorem \ref{3121} $x_1\ge y_1+y_2\ge 3$.
 Then $G$ contains the difference graph with vertex-eigenvector $(3,\sum_{i=1}^{k-1}x_i-3,x_k;1,2,\sum_{j=1}^{k}y_j-3)$ as a proper subgraph. Thus  $$\left \{ \begin{array}{ll}
m\ge 3(\sum_{i=1}^{k-1}x_i+\sum_{j=1}^{k}y_j)-9+x_k;\\
n=\sum_{i=1}^{k}x_i+\sum_{j=1}^{k}y_j, \end{array}\right.$$
which implies that
$x_k\ge \frac{3n-m-9}{2}\ge \frac{3n-m-8}{2}$ as $3n-m-9$ is odd. Moreover, recall that $y_2\ge 2$, then $m-\sum_{i=1}^{k}x_i\ge m-n+4$. Consequently,
\begin{equation*}
\begin{array}{lll}\mathbf{a}_4(G)&\ge & \frac{(3n-m-8)(m-n+4)}{2}+(1+y_2)x_{k-1}y_kx_1\\
&\ge & \frac{(3n-m-8)(m-n+4)}{2}+9\\
&>&\frac{(3n-m-10)(m-n+3)}{2}+m-n+2,
\end{array}
\end{equation*}
 which  contradicts to that $G$ has minimal $4$-Sachs number in $\mathfrak{B}_{n,m}$. Consequently, $y_2=1$.

\noindent {\bf Assertion 3.} $x_{k-1}=1$.

 By {\bf Assertion 2} the vertex-eigenvector of $G$ is $(x_1,x_2,\ldots,x_{k-1},x_k;1,1,\ldots,y_{k-1},y_k)$. We divide two steps to prove that $x_{k-1}=1$.
  Firstly, we show that $x_{k-1}\le y_3$. Assume to the contrary that $x_{k-1}> y_3$. Let $x_{k-1}=p(1+y_3)+q$ with $p\ge 1$ and $0\le q\le y_3$. Note that $(\sum_{i=1}^{k-1}x_i, 2)$ is a outer corner and $(\sum_{i=1}^{k-2}x_i+1,3)$ is an inner corner, then by Corollary \ref{316}  $$2\sum_{i=1}^{k-1}x_i\ge 3(\sum_{i=1}^{k-2}x_i+1),$$ that is, $$2x_{k-1}\ge \sum_{i=1}^{k-2}x_i+3.\eqno{(5.4)}$$ Let now $G'$ be the the difference graph with vertex eigenvector $(x_1,x_2,\ldots,x_{k-2}+p, x_{k-1}-py_3-p,x_k+py_3;1,1,y_3,\ldots,y_{k-1},y_k)$. Then the
   characteristic matrix of $G'$ is
$$T(G')=\left(
                        \begin{array}{ccccc}
                         x_1& x_1  &x_1y_3&\ldots & x_1y_k \\
                         \ldots& \ldots  &\ldots&\ldots & \ldots \\
        x_{k-2}+p& x_{k-2}+p &(x_{k-2}+p)y_3&\ldots & 0 \\
        x_{k-1}-py_3-p& x_{k-1}-py_3-p &0&\ldots & 0 \\
        x_k+py_3& 0 &0&\ldots & 0
                        \end{array}
                      \right).$$
One find that $G$ and $G'$ have the same number of edges. By a directly calculation, we have
Thus \begin{equation*}
\begin{array}{lll}\mathbf{a}_4(G')-\mathbf{a}_4(G)&=&(\sum_{i=1}^{k-2}x_i+p)y_3(2x_{k-1}+x_k-py_3-2p)-(\sum_{i=1}^{k-1}x_i)x_k \\
&&-(\sum_{i=1}^{k-2}x_i)y_3(2x_{k-1}+x_k)+(\sum_{i=1}^{k-1}x_i-py_3)(x_k+py_3)\\
&=& py_3[2x_{k-1}-2py_3-2p-(y_3+2)(\sum_{i=1}^{k-2}x_i)+\sum_{i=1}^{k-1}x_i]\\
&=& py_3[2q+x_{k-1}-(y_3+1)(\sum_{i=1}^{k-2}x_i)]\\
&\le&py_3[2q+x_{k-1}+(y_3+1)(3-2x_{k-1})]~~~~~~~~~~(by~~ (5.4))\\
&\le&py_3(2q-y_3x_{k-1})~~(as~x_{k-1}\ge 3 ~by~(5.4)~ and~ Theorem ~\ref{3121})\\
&=&q(1-y_3)+q-py_3-pqy_3\\
&<&0,
\end{array}
\end{equation*}which implies that $x_{k-1}\le y_3$. \vspace{2mm}

Assume now that $2\le x_{k-1}\le y_3$. Then $y_3\ge 2$. Let $\{P_j(\sum_{i=1}^{k-1}x_i-j,2)|j=0,\ldots,x_{k-1}-1\}$ and $\{Q_j(\sum_{i=1}^{k-2}x_i+1,3+j)|j=0,\ldots,x_{k-1}-1\}$. Let also $G_0=G$ and $G_j=G_{j-1}-P_j+Q_j$ for $j=1,\ldots,x_{k-1}-1$. Then one can find that for each $j$ $P_j$ is a outer corner and $Q_j$ is an inner corner of $G_{j-1}.$ Let $x:=\sum_{i=1}^{k-2}x_i$. Then
\begin{equation*}
\begin{array}{lll}&&2\sum_{j=0}^{x_{k-1}-1}(\sum_{i=1}^{k-1}x_i-j)-(\sum_{i=1}^{k-2}x_i+1)\sum_{j=0}^{x_{k-1}-1}(3+j)\\
&=&x_{k-1}[(2x+1)-\frac{(5+x_{k-1})(x+1)}{2}]\\
&<&0.
\end{array}
\end{equation*}
Then $\mathbf{a}_4(G_{x_{k-1}})<\mathbf{a}_4(G)$ by Corollary \ref{317}, a contradiction. Consequently, $x_{k-1}=1$.

\noindent {\bf Assertion 4.} The character of $G$ is $3$.

 Assume to the contrary that the character $k\ge 4$. Set $y=\sum_{i=1}^{k-3}x_i.$
Note that  $(y+x_{k-2}+1, 2)$ is a outer corner and $(y+1,3+y_3)$ is an inner corner, then from Lemma \ref{316} $$2(y+x_{k-2}+1)\ge (3+y_3)(y+1),$$
that is $$x_{k-2}\ge (1+y)(1+y_3)>yy_3.\eqno{(5.5)}$$
On the other hand, note that $(y+x_{k-2}+1, 3)$ is an inner corner and $(y,3+y_3)$ is a outer corner, applying Lemma \ref{316} again
a contradiction to (5.5) is yielded. Thus the character of $G$ is $3$.

\noindent {\bf Assertion 5.} $y_3=1$.

From Assertions 1 to 4 the vertex-eigenvector of $G$ is $(x_1,1,x_3;1,1,y_3)$. Then  $$\left \{ \begin{array}{ll}
m=  x_1(2+y_3)+2+x_3;\\
n=3+y_3+x_1+x_3. \end{array}\right.$$
From which we have $$\left \{ \begin{array}{ll}
x_3=n-4-y_3-\frac{m-n}{1+y_3}\ge \frac{3n-m-10}{2};\\
m-x_1-1-x_3\ge m-n+3;\\
x_1y_3=(\frac{m-n}{1+y_3}+1)y_3 \end{array}\right.$$
 with equality  if and only if $y_3=1$. Thus
 $$\mathbf{a}_4(G)=x_3(m-x_1-x_2-x_3)+2x_1y_3\ge \frac{(3n-m-10)(m-n+3)}{2}+m-n+2$$with equality  if and only if $y_3=1$.\vspace{3mm}

Up to now, we show that the difference graph $G$ having minimal $4$-Sachs number has vertex-eigenvector $w_3=(\frac{m-n+2}{2},1,\frac{3n-m-10}{2};1,1,1)$ if $y_1=1$ and $3n-m-9$ is odd.

Comparing $\mathbf{a}_4(G,w_i)$($i=1,2,3,4$), we have $\mathbf{a}_4(G,w_1)<min\{\mathbf{a}_4(G,w_2),\mathbf{ a}_4(G,w_3)\}$ if $m< \frac{7n}{3}-7$; $\mathbf{a}_4(G,w_1)=\mathbf{a}_4(G,w_2)$ if $m= \frac{7n}{3}-7$;
$\mathbf{a}_4(G,w_1)>\mathbf{a}_4(G,w_2)$ if $m> \frac{7n}{3}-7$ and $3n-m-9$ is even; and $\mathbf{a}_4(G,w_1)>\mathbf{a}_4(G,w_3)$ if $m> \frac{7n}{3}-7$ and $3n-m-9$ is odd.
Consequently, the proof is complete. \hfill $\blacksquare$

Combining Theorem \ref{319} and Theorem \ref{320}.
\begin{theorem} \label{320}Let $n>6$ and $n-1\le m\le  3(n-3)$. Then
 $$\bar{\mathbf{a}}_4(\mathfrak{B}_{n,m})=\left \{ \begin{array}{ll}
 (2n-4-m)(m-n+1), & {\rm if  \mbox{ } n-1\le m\le  2(n-2) };\\
2(3n-9-m)(m-2n+6), & {\rm if  \mbox{ } m< \frac{7n}{3}-7 };\\
\frac{(3n-9-m)(m-n+3)}{2}, & {\rm if  \mbox{ } m> \frac{7n}{3}-7 \mbox{ } and \mbox{ }3n-m-9\mbox{ }is \mbox{ } even};\\
\frac{(3n-m-10)(m-n+3)}{2}+m-n+2, & {\rm if  \mbox{ } m> \frac{7n}{3}-7 \mbox{ } and \mbox{ }3n-m-9\mbox{ }is \mbox{ } odd.} \end{array}\right.$$
\end{theorem}

\section{Another formula for  $4$-Sachs number of difference graphs}
In this section, we first establish another formula for  $4$-Sachs number of difference graphs. Then we express the problem of computing the minimal $4$-Sachs number  as a combinatorial optimization problem, which relates to the partitions of positive integers.
\begin{theorem}\label{331} Let $G$ be a difference graph with  Young matrix $Y(G)$, defined as above.  Denote by $r_1,r_2,\ldots, r_h$ the row sum of $Y$, respectively. Then
$$\mathbf{a}_4(G)=\sum_{i=1}^{h-1}\sum_{j=i+1}^h(r_i- r_{i+1})i r_j.$$
\end{theorem}
\noindent {\bf Proof.} Suppose that $(X_1,X_2,\ldots,X_k; Y_1,Y_2,\ldots,Y_k)$ is the vertex-bipartition of $G$. If $r_i- r_{i+1}>0$ for some $i$, then there exists an integer $t$ such that $\sum_{s=1}^t|X_s|=i$ thus $(c_i-c_{i+1})i$ denotes the number of edges contained in $G[\cup_{j=1}^tX_j;Y_{k-t+1}]$ and $\sum_{j=i+1}^h r_j$ denotes the number of edges contained in the difference graph $G\overline{[\cup_{j=1}^tX_j;Y_{k-t+1}]}$. Then combining with Theorem \ref{32}, the result follows.\hfill $\blacksquare$\vspace{2mm}

A partition of a positive integer $n$ is any non-increasing sequence of positive integers whose sum is $n$. The problem on partitions of  positive integers was first studied by G. W. Leibniz; see  \cite{an,ge}.  Let $(x_1,x_2,\ldots,x_k; y_1,y_2,\ldots,y_k)$ be the vertex-eigenvector of the difference graph $G$  with  $ \sum_{i=1}^kx_i\ge \sum_{i=1}^ky_i$. As we known that the row sum sequence $r_1,r_2,\ldots, r_h$ of its Young diagram is a non-increase positive sequence satisfying $r_1+h=n-1$
and $\sum_{i=1}^hr_i=m$. Then the problem of computing the minimal $4$-Sachs number, in $\mathfrak{B}_{n,m}$, can be expressed as a optimization problem related to the partition of a positive integer with the following restrictive conditions.

$$\min \sum_{i=1}^{h-1}\sum_{j=i+1}^h(r_i- r_{i+1})i r_j.$$
s.t.
$$\left \{ \begin{array}{l}
r_1\ge r_2\ge \ldots\ge r_h>0;\\
r_1+h=n-1;\\
h\ge \lfloor\frac{n-1}{2}\rfloor;\\
\sum_{i=1}^hr_i=m.
 \end{array}\right.$$
Unfortunately, we  say nothing on minimal $4$-Sachs number from the restrictive conditions above.


\baselineskip=0.21in

\begin{thebibliography}{99}
\bibitem{an} G. Andrews, The Theory of Partitions. \emph{Addison-Wesley Publishing Company},  1976.
\bibitem{as}Y. Ashkenazi, $C_3$ saturated graphs, \emph{ Discrete Mathematics},  {\bf 297}( 2005)  152-158.

\bibitem{big} N. L. Biggs, Algebraic Graph Theory. \emph{Cambridge University Press,} Cambridge, 1993.
\bibitem{br}R. A. Brualdi, H. J. Ryser. Combinatorial Matrix Theory. \emph{Cambridge University Press}, Cambridge,1991.

\bibitem{be}B. Bollob$\acute{a}$s, P. Erd$\ddot{o}$s, On a Ramsey-Tur$\acute{a}$n type problem, \emph{J. Combin. Theory B} {\bf 21}(1976) 166-168.




  \bibitem{cds}  D. Cvetkovi$\acute{c}$, M. Doob, H. Sachs, Spectra of Graphs, \emph{Academic Press}, New York, 1980.

\bibitem{crs} D. Cvetkovi$\acute{c}$, P. Rowlinson, S. Simi$\acute{c}$, An introduction to the theory of graph spectra, \emph{Cambridge University Press}, Cambridge, 2009.

\bibitem{ge}I.M. Gessel, Counting paths in Young's lattice, \emph{ Journal Statistical Planning and Inference} {\bf 34}(1993) 125-134.


\bibitem{gr} C. Godsil, G. Royle, Algebraic Graph Theory, \emph{Springer-Verlag,} New York Inc.; 2001.


\bibitem{g} S. C. Gong, On the rank of a real skew symmetric matrix described by an oriented graph, \emph{ Linear and Multilinear Algebra,} {\bf 65}(2017) 1934-194.

%
%

\bibitem{hps}P. L. Hammer, U. N. Peled,  X. R. Sun, Difference graphs, \emph{ Discrete Appl. Math.} {\bf 28}(1990) 35-44.

%


\bibitem{k}L. K$\acute{a}$szonyi, Z. Tuza, Saturated graphs with minimal number of edges, \emph{ J. Graph Theory} {\bf 10}(1986) 203-210.

\bibitem{lmt}J. Lazzarin, O. F. M$\acute{a}$rquez, F. C. Tura, No threshold graphs are cospectral, \emph{ Linear Algebra Appl.} {\bf 560}(2019) 133-145.


\bibitem{lp} L. Lov$\acute{a}$sz, M. Plummer,\emph{ Matching Theory}, Ann. Discrete Math., vol.
29, \emph{North-Holland}, New York, 1986.


\bibitem{mp}N. V. R. Mahadev, U. N. Peled. Threshold Graphs and Related Topics. Elsevier Publishers, 1995.
\bibitem{m}A. Mowshowitz, The characteristic polynomial of a graph, \emph{ J. Combin. Theory Ser. B} {\bf 12}(1972) 177-193.

\bibitem{sw} A. J. Schwenk, R. J. Wilson. On the eigenvalues of a graph. Selected Topics in Graph Theory,  \emph{Academic Press,} New York, 1978.



\bibitem{kr} L. Keough, A.J. Radcliffe, Graphs with the fewest matchings, \emph{ Combinatorica,}
 {\bf 36(6)}(2016), 703-723.

 \bibitem{y} M. Yannakakis, The complexity of the partial order dimension problem, \emph{ SIAM J. Algebraic Discrete
Methods} {\bf 3}(1982) 351-358.
\end{thebibliography}
\end{document}